\documentclass[12pt, oneside, leqno]{amsart}

\usepackage{amsthm}

\newtheorem{thm}{Theorem}[section]    
\newtheorem{prop}[thm]{Proposition}

\theoremstyle{definition}
\newtheorem{defn}[thm]{Definition}    
\newtheorem{rem}[thm]{Remark}             
\newtheorem{exa}[thm]{Example}

\theoremstyle{mytheoremstyle}

\newtheorem{cor}[thm]{Corollary}

\newcommand*{\longhookrightarrow}{\ensuremath{\lhook\joinrel\relbar\joinrel\rightarrow}}

\usepackage{amssymb}
\usepackage{amsmath}
\usepackage{extpfeil}
\usepackage{tikz}
\usepackage{tikz-cd}
\usepackage{here}
\usepackage[figuresright]{rotating}

\def\qbb{{\mathbb{Q}}}
\def\rbb{{\mathbb{R}}}

\def\ccal{{\mathcal{C}}}
\def\dcal{{\mathcal{D}}}

\def\vcal{{\mathcal{V}}}
\def\wcal{{\mathcal{W}}}

\def\xfra{{\mathfrak{X}}}

\def\czero{{\mathcal{C}^0}}

\numberwithin{equation}{section}

\newextarrow{\xbigtoto}{{20}{20}{20}{20}}
{\bigRelbar\bigRelbar{\bigtwoarrowsleft\rightarrow\rightarrow}}
\usepackage[top=30truemm,bottom=30truemm,left=30truemm,right=30truemm]{geometry}

\begin{document}

\title{$C^\infty$-manifolds with skeletal diffeology}
\author{Hiroshi Kihara}

\maketitle
\renewcommand{\thefootnote}{}

\footnote{Center for Mathematical Sciences, University of Aizu, 
	Tsuruga, Ikki-machi, Aizu-Wakamatsu City, Fukushima, 965-8580, Japan} \footnote{\textbf{e-mail: kihara@u-aizu.ac.jp}}

\renewcommand{\thefootnote}{\arabic{footnote}}
\setcounter{footnote}{0}
\vspace{-5ex}
\begin{abstract}
We formulate and study the notion of $d$-skeletal diffeology, which generalizes that of wire diffeology, introducing the dual notion of $d$-coskeletal diffeology. We first show that paracompact finite-dimensional $C^\infty$-manifolds $M_d$ with $d$-skeletal diffeology inherit good topologies and smooth paracompactness from $M$. We then study the pathology of $M_d$. Above all, we prove the following: For $d<{\rm dim}\ M$, every immersion $f:M\longrightarrow N$ is isolated in the diffeological space $\dcal(M_d, N_d)$ of smooth maps and the $d$-dimensional smooth homotopy group of $M_d$ is uncountable.
\end{abstract}

\section{Introduction}\label{intro}

Let $X$ be a diffeological space. Then, we can consider the new diffeology on the set $X$ generated by the $1$-dimensional plots (see \cite[1.66]{IZ}); this diffeology, introduced by Souriau, is called the wire diffeology (or spaghetti diffeology). Similarly, we can consider the new diffeologies on the set $X$ generated by the plots of dimension $\le d$ for $d=2,3,\ldots$ (see \cite[1.10]{IZ}).

$C^{\infty}$-manifolds with wire diffeology are well-known pathological diffeological spaces; the definition and the feeling that they are pathological are shared among experts. However, no general theory and few pathological phenomena have been reported. Thus, we develop the general theory of wire diffeology and its higher dimensional analogues, introducing their dual notions, and report new pathological phenomena. The main results are outlined below.

We formulate the notions of wire diffeology and its higher dimensional analogues as the $d$-skeleton functors
\[
\cdot_d : \dcal \longrightarrow \dcal
\]
$(0<d<\infty)$, where $\dcal$ denotes the category of diffeological spaces. We also construct the $d$-coskeleton functors
\[
\cdot_{(d)}: \dcal \longrightarrow \dcal,
\]
showing that $\cdot_d : \dcal \rightleftarrows \dcal: \cdot_{(d)}$ is an adjoint pair for any $d$ (Proposition \ref{adjoint}). This is a diffeological analogue of the adjoint pair of the $d$-skeleton and $d$-coskeleton functors for simplicial sets. (See Section 3 for details.)

First, we state that the $d$-skeleton functor preserves underlying topologies and some smooth topological properties of diffeological spaces. The basics of the underlying topological space functor $\widetilde{\cdot}: \dcal \longrightarrow \czero$ and those of $\dcal$-paracompactness and $\dcal$-regularity are presented in Remark \ref{convenrem}(2) and Section 2.3, respectively. Throughout this paper, $C^\infty$-manifolds are ones in the sense of \cite[Section 27]{KM}; see Section 2.2. (However, if the reader is unfamiliar with convenient calculus, then he/she may assume that $C^\infty$-manifolds are paracompact and finite-dimensional.)

\begin{prop}\label{k-trun}
	Let $X$ be a diffeological space, and $d$ an integer with $0<d<\infty$.
	\begin{itemize}
		\item[{\rm (1)}] The map $id:X_{d}\longrightarrow X$ is a smooth map such that $\widetilde{id}:\widetilde{X_{d}}\longrightarrow\widetilde{X}$ is a homeomorphism.
		\item[{\rm (2)}] The map $id:X_d \longrightarrow X$ induces the diffeomorphism 
		\[
		\dcal(X_d,N) \overset{id^\sharp}{\longleftarrow} \dcal(X,N) 
		\]
		for any $C^\infty$-manifold $N$.
		\item[{\rm (3)}]  $X$ is $\mathcal{D}$-$paracompact$ (resp. $\dcal$-regular) if and only if $X_d$ is $\dcal$-paracompact (resp. $\dcal$-regular).
	\end{itemize}
\end{prop}

The obvious analogue of Proposition \ref{k-trun} holds for $X_{(d)}$ and the map $id: X \longrightarrow X_{(d)}$ (see Remark \ref{(co)unit}). However, since $M=M_{(d)}$ for any $C^\infty$-manifold $M$ and any $d>0$ (Proposition \ref{mfdcosk}), we are mainly concerned with (paracompact finite-dimensional) $C^\infty$-manifolds $M_d$ with skeletal diffeology.

As observed below, $C^\infty$-manifolds with skeletal diffeology are certainly pathological; see \cite[Example 3.22(1)]{CW16a} and \cite[Section 6]{CW20} for the pathology of wire diffeology. On the other hand, such diffeological spaces have good underlying topologies, along with good smooth topological properties such as $\dcal$-paracompactness, presenting a sharp contrast to other pathological diffeological spaces such as irrational tori. Thus, $C^\infty$-manifolds with skeletal diffeology are expected to be useful as counterexamples to precise conjectures in diffeology which assume good topological properties or good smooth topological properties (see Remark \ref{use}).

We give main theorems, describing smooth homotopical pathology of $C^\infty$-manifolds with skeletal diffeology. In particular, the following theorem and its corollaries show that smooth maps between $C^{\infty}$-manifolds behave well for the skeleton functor $\cdot_{d}:\mathcal{D}\longrightarrow\mathcal{D}$, but that smooth homotopies between them don't.

A smooth map $f:M \longrightarrow N$ between $C^\infty$-manifolds is called an {\sl immersion} if $T_x f : T_x M \longrightarrow T_{f(x)} N$ is injective for any $x\in M$ (see \cite[p. 284]{KM}).

\begin{thm}\label{immersion}
	Let $M$ and $N$ be $C^{\infty}$-manifolds, and $d$ an integer with $0<d<\infty$.
	\begin{itemize}
		\item[{\rm (1)}] The equalities
		\[
		C^{\infty}(M,N)=\mathcal{D}(M,N)=\mathcal{D}(M_{d},N_{d}) 
		\]
		hold in Set.
		\item[{\rm (2)}] If $f:M\longrightarrow N$ is an immersion, then the set $\{f\}$ is a connected component of $\mathcal{D}(M_{d},N_{d})\ \rm{for}\ \it{d}<\rm{dim}\ \it{M}$.
	\end{itemize}
\end{thm}

For a diffeological space $F$, the automorphism group ${\rm Aut}_\dcal(F)$ is just the diffeomorphism group of $F$. Endow the group ${\rm Aut}_\dcal(F)$ with the initial structure for the two maps
\[
{\rm Aut}_\dcal (F) \xbigtoto[inv]{incl} \dcal(F,F),
\]
where $incl(\varphi) = \varphi$ and $inv(\varphi)=\varphi^{-1}$. Then, we can easily see that ${\rm Aut}_\dcal (F)$ is a diffeological group. (See \cite[Section 5.3]{smh} and \cite[1.61]{IZ}).

The following result is derived from Theorem \ref{immersion}.

\begin{cor}\label{auto}
	Let $M$ be a $C^{\infty}$-manifold, and $d$ an integer with $0<d <{\rm dim}\ M$. Then, ${\rm Aut}_{\mathcal{D}}(M_{d})$ is just the diffeomorphism group of $M$, endowed with the discrete diffeology.
\end{cor}

We can also derive the following result from Theorem \ref{immersion}.

\begin{cor}\label{homotopydiffeo}
	Let $M$ and $N$ be $C^\infty$-manifolds and let $0<d<{\rm min}({\rm dim}\ M, {\rm dim}\ N)$. Then, the following are equivalent:
	\begin{itemize}
		\item[{\rm (i)}] $M_d$ is smoothly homotopy equivalent to $N_d$.
		\item[{\rm (ii)}] $M_d$ is diffeomorphic to $N_d$.
		\item[{\rm (iii)}] $M$ is diffeomorphic to $N$.
	\end{itemize}
\end{cor}

Every $C^\infty$-manifold is locally contractible as a diffeological space \cite[Lemma 4.17]{KM}, and hence as a topological space (see \cite[Section 2.4]{origin}). However, this is not the case for $C^\infty$-manifolds with skeletal diffeology.

\begin{cor}\label{loccontr}
	Let $M$ be a $C^\infty$-manifold and $d$ an integer with $0<d<{\rm dim}\ M$. Then, the $C^\infty$-manifold $M_d$ with $d$-skeletal diffeology is a non-locally contractible diffeological space whose underlying topological space is locally contractible.
\end{cor}

Recall the natural homomorphism
\[
\pi^\dcal_i (X, x_0) \longrightarrow \pi_i (\widetilde{X},x_0)
\]
between the smooth and topological homotopy groups from \cite[Section 3]{CW}. Note that $\widetilde{X}_d = \widetilde{X}$ (Proposition \ref{k-trun}(1)). See Section 2.2 and Remark \ref{Cparacpt} for the notions of $c^\infty$-topology and convenient vector space and that of $C^\infty$-regularity, respectively.

\begin{thm}\label{notVD}
	Let $M$ be a $C^\infty$-manifold and let $0<d< {\rm dim}\ M$ and $x_0 \in M$. 
	\begin{itemize}
		\item[{\rm (1)}] The natural homomorphism 
		\[
		\pi^\dcal_i (M_d, x_0) \longrightarrow \pi_i(\widetilde{M}, x_0)
		\]
		is an isomorphism for $i<d$ and an epimorphism for $i=d$.
		\item[{\rm (2)}] Suppose that $M$ satisfies one of the following conditions:
		\begin{itemize}
			\item[{\rm (i)}] $M$ is a $c^\infty$-open set of a convenient vector space.
			\item[{\rm (ii)}] $M$ is $C^\infty$-regular.
		\end{itemize}
		Then, the kernel of the epimorphism 
		\[
		\pi^\dcal_d (M_d, x_0) \longrightarrow \pi_d (\widetilde{M}, x_0)
		\]
		is uncountable.
	\end{itemize}
\end{thm}

We end this section with a little more detail on why we study $C^{\infty}$-manifolds with skeletal diffeology.

\begin{rem}\label{use}
In diffeology, one tries to generalize the notions and results for paracompact finite-dimensional $C^{\infty}$-manifolds. Then, typical pathological diffeological spaces such as irrational tori $T_{\Gamma}$ and $\rbb / \qbb$ are often useful as counterexamples to naive conjectures. However, paracompact finite-dimensional $C^{\infty}$-manifolds have good topological properties such as local contractibility and paracompactness, and good smooth topological properties such as (smooth) local contractibility and smooth paracompactness. Thus, in various situations, one should make precise conjectures assuming good topological properties or good smooth topological properties. But the above mentioned typical pathological diffeological spaces cannot be counterexamples to such precise conjectures because of their too bad underlying topologies.

Thus, we need pathological diffeological spaces with good underlying topologies and some good smooth topological properties. Hence, we should pay attention to $C^{\infty}$-manifolds $M_{d}$ with skeletal diffeology and diffeological spaces of the form $RZ$ (see Remark \ref{convenrem}(2)). For simplicity, let $M$ be a paracompact finite-dimensional $C^{\infty}$-manifold. Then, the following hold:
\begin{itemize}
\item For $0 < d < \mathrm{dim} \ M$, $M_{d}$ is locally contractible and paracompact as a topological space, and is smoothly paracompact but non-locally contractible as a diffeological space (see Proposition \ref{k-trun} and Corollary \ref{loccontr}).
\item $R\widetilde{M}$ is locally contractible and paracompact as a topological space, and is locally contractible but non-smoothly paracompact as a diffeological space (see \cite[Proposition 2.10]{origin} and \cite[Step 2 of the proof of Corollary 1.6(2) and Remark A.8]{smh}).
\end{itemize}

\noindent Thus, $C^{\infty}$-manifolds $M_{d}$ with skeletal diffeology, along with diffeological spaces of the form $RZ$ are expected to be useful as counterexamples to precise conjectures which assume good topological properties or good smooth topological properties (see our forthcoming paper concerning de Rham caluculus on diffeological spaces \cite{deRham}).  
\end{rem}

\section{Diffeological spaces and $C^\infty$-manifolds}
We recall the convenient properties of diffeological spaces and arc-generated spaces, along with the fact that the category $C^\infty$ of (separated) $C^\infty$-manifolds can be fully faithfully embedded into the category $\dcal$ of diffeological spaces. We also make a brief review on smooth paracompactness and smooth regularity.

\subsection{Categories $\dcal$ and $\czero$}
In this subsection, we summarize the convenient properties of the category $\dcal$ of diffeological spaces, recalling the adjoint pair $\widetilde{\cdot}: \dcal \rightleftarrows \czero: R$ of the underlying topological space functor and its right adjoint; see \cite{IZ} and \cite{origin} for full details.
\par\indent
Let us begin with the definition of a diffeological space. A {\sl parametrization} of a set $X$ is a (set-theoretic) map $p: U \longrightarrow X$, where $U$ is an open subset of $\rbb^{n}$ for some $n$.

\begin{defn}\label{diffeological}
	\begin{itemize}
		\item[(1)] A {\sl diffeological space} is a set $X$ together with a specified set $D_X$ of parametrizations of $X$ satisfying the following conditions:
		\begin{itemize}
			\item[(i)](Covering)  Every constant parametrization $p:U\longrightarrow X$ is in $D_X$.
			\item[(ii)](Locality) Let $p :U\longrightarrow X$ be a parametrization such that there exists an open cover $\{U_i\}$ of $U$ satisfying $p|_{U_i}\in D_X$. Then, $p$ is in $D_X$.
			\item[(iii)](Smooth compatibility) Let $p:U\longrightarrow X$ be in $D_X$. Then, for every $n \geq 0$, every open set $V$ of $\rbb^{n}$ and every smooth map $F  :V\longrightarrow U$, $p\circ F$ is in $D_X$.
		\end{itemize}
		The set $D_X$ is called the {\sl diffeology} of $X$, and its elements are called {\sl plots}.
		\item[(2)] Let $X=(X,D_X)$ and $Y=(Y,D_Y)$ be diffeological spaces, and let $f  :X\longrightarrow Y$ be a (set-theoretic) map. We say that $f$ is {\sl smooth} if for any $p\in D_X$, \ $f\circ p\in D_Y$.
	\end{itemize}
\end{defn}
The convenient properties of $\dcal$ are summarized in the following proposition. Recall that a topological space $X$ is called {\sl arc-generated} if its topology is final for the continuous curves from $\rbb$ to $X$, and let $\czero$ denote the category of arc-generated spaces and continuous maps. See \cite[pp. 230-233]{FK} for initial and final structures with respect to the underlying set functor.

\begin{prop}\label{conven}
	\begin{itemize}
		\item[$(1)$] The category ${\dcal}$ has initial and final structures with respect to the underlying set functor. In particular, ${\dcal}$ is complete and cocomplete.
		\item[$(2)$] The category $\mathcal{D}$ is cartesian closed.
		\item[$(3)$] The underlying set functor $\dcal \longrightarrow Set$
		is factored as the underlying topological space functor
		$\widetilde{\cdot}:\dcal \longrightarrow \czero$
		followed by the underlying set functor
		$\czero \longrightarrow Set$.
		Further, the functor
		$\widetilde{\cdot}:\dcal \longrightarrow \czero$
		has a right adjoint
		$R:\czero \longrightarrow \dcal$.
	\end{itemize}
	\begin{proof}
		See \cite[p. 90]{CSW}, \cite[pp. 35-36]{IZ}, and \cite[Propositions 2.1 and 2.10]{origin}.
	\end{proof}
\end{prop}
The following remark relates to Proposition \ref{conven}.
\begin{rem}\label{convenrem}
	{\rm (1)} Let $\xfra$ be a concrete category (i.e., a category equipped with a faithful functor to $Set$); the faithful functor $\xfra \longrightarrow Set$ is called the underlying set functor. See \cite[Section 8.8]{FK} for the notions of an $\xfra$-embedding, an $\xfra$-subspace, an $\xfra$-quotient map, and an $\xfra$-quotient space. $\dcal$-subspaces and $\dcal$-quotient spaces are usually called diffeological subspaces and quotient diffeological spaces, respectively.
	
	{\rm (2)} For Proposition \ref{conven}(3), recall that the underlying topological space $\tilde{A}$ of a diffeological space $A = (A, D_A)$ is defined to be the set $A$ endowed with the final topology for $D_A$ and that $R$ assigns to an arc-generated space $Z$ the set $Z$ endowed with the diffeology 
	$$
	D_{RZ} = \text{\{continuous parametrizations of } Z \text{\}.}
	$$
	Then, we can easily see that $\tilde{\cdot} \circ R = Id_{\ccal^0}$ and that the unit $A \longrightarrow R\widetilde{A}$ of the adjoint pair $(\widetilde{\cdot}, R)$ is set-theoretically the identity map.
	
	{\rm (3)} The notion of an arc-generated space is equivalent to that of a $\Delta$-generated space (see \cite[Section 3.2]{CSW}, \cite[Section 2.2]{origin}). The categories $\dcal$ and $\czero$ share convenient properties (1) and (2) in Proposition \ref{conven}, which often enables us to deal with $\dcal$ and $\czero$ simultaneously (see \cite{smh}). See \cite[Remark 2.4]{smh} for the reason why $\czero$ is the most suitable category as a target category of the underlying topological space functor for diffeological spaces.
\end{rem}

\subsection{$C^\infty$-manifolds as diffeological spaces}
In convenient calculus \cite{KM}, $C^\infty$-manifolds are defined by gluing $c^\infty$-open sets of convenient vector spaces. Here, convenient vector spaces are locally convex spaces satisfying a weak completeness condition and the $c^\infty$-topology of a locally convex space is defined to be the final topology for the smooth curves. The $c^\infty$-topology of a locally convex space $E$ is finer than the original locally convex topology, but the two topologies coincide if $E$ is metrizable. Of course, finite-dimensional $C^\infty$-manifolds in convenient calculus are just finite-dimensional $C^\infty$-manifolds in the usual sense. (See \cite[the beginning of Section 2.2]{smh} for a brief explanation of convenient calculus. See also \cite{KM} for full details.)

Throughout this paper, $C^{\infty}$-manifolds are ones in convenient calculus (i.e., $C^\infty$-manifolds in the sense of \cite[Section 27]{KM}); we assume that $C^{\infty}$-manifolds are separated as arc-generated spaces. Separatedness is stronger than $T_1$-axiom and weaker than Hausdorff property. (See \cite[the end of Section 2.1]{smh} for the notion of separatedness.) For finite-dimensional $C^\infty$-manifolds, separatedness is equivalent to Hausdorff property (see \cite[Section 2.2]{origin}).

The category $C^\infty$ of $C^\infty$-manifolds can be regarded as a full subcategory of the category $\dcal$ of diffeological spaces. In fact, the fully faithful embedding
\[
I: C^\infty \longhookrightarrow \dcal 
\]
is defined to assign to a $C^\infty$-manifold $M$ the set $M$ endowed with the diffeology 
\[
D_{IM} = \text{\{$C^\infty$-parametrizations of $M$\}}
\]
(see \cite[Section 2.2]{smh}). The diffeological space $IM$ is often denoted simply by $M$ if there is no confusion in context.

\subsection{Smooth paracompactness and smooth regularity}
We recall the notions of smooth paracompactness and smooth regularity from \cite[Section 5.1]{smh} and \cite[Sections 14 and 16]{KM}.

\begin{defn}\label{Dparacpt}
	Let $X$ be a diffeological space.
	\begin{itemize}
		\item[{\rm (1)}] $X$ is called $\dcal$-paracompact if for any open covering $\{U_i\}$ of $X$, there exists a smooth partition of unity $\{\varphi_i : X\longrightarrow \rbb\}$ subordinate to $\{U_i\}$ (i.e., ${\rm supp}\ \varphi_i \subset U_i$).
		\item[{\rm (2)}] $X$ is called $\dcal$-regular if for any $x\in X$ and any open neighborhood $U$ of $x$, there exists a smooth function $f: X\longrightarrow \rbb$ such that $f(x)=1$ and ${\rm carr}(f) \subset U$, where 
		\[
		{\rm carr}(f) = \{z\in X\ |\ f(z) \neq 0 \}. 
		\]
	\end{itemize}
	In the definition of $\dcal$-regularity, ${\rm carr} (f)$ can be replaced with ${\rm supp}(f)$. If $X$ is a $\dcal$-paracompact diffeological space with $\widetilde{X}$ $T_1$-space, then $X$ is $\dcal$-regular. We can easily see that separation axiom $T_1$ is equivalent to Hausdorff condition under $\dcal$-regularity (hence under $\dcal$-paracompactness).
\end{defn}

\begin{rem}\label{Cparacpt}
	Similarly, we can define the notions of $C^\infty$-paracompactness and $C^\infty$-regularity for $C^\infty$-manifold. For a $C^\infty$-manifold $M$, $C^\infty$-paracompactness and $\dcal$-paracompactness (resp. $C^\infty$-regularity and $\dcal$-regularity) are equivalent, so these notions are often referred to as smooth paracompactness (resp. smooth regularity).
	
	Since $C^\infty$-manifolds are assumed to be separated as arc-generated spaces, $C^\infty$-regular $C^\infty$-manifolds, and hence $C^\infty$-paracompact $C^\infty$-manifolds are Hausdorff.
	
	Paracompact finite-dimensional $C^\infty$-manifolds are $C^\infty$-paracompact, and hence $C^\infty$-regular (see \cite[Theorem 7.3]{BJ}). Further, such $C^\infty$-manifolds are hereditarily $C^\infty$-paracompact by \cite[Corollary 11.17]{smh} or \cite[Theorem 2]{Gauld}. In \cite[Section 11.4]{smh}, it is shown that many important $C^\infty$-manifolds are (hereditarily) $C^\infty$-paracompact.
\end{rem}

\section{Skeletal and coskeletal diffeologies}
In this section, we introduce the $d$-skeleton functor $\cdot_d:\dcal \longrightarrow \dcal$, along with the $d$-coskeleton functor $\cdot_{(d)}: \dcal \longrightarrow \dcal$, and prove Proposition \ref{k-trun}.


First, we define the $d$-skeleton functor
\[
\cdot_d: \mathcal{D} \longrightarrow \mathcal{D}
\]
to assign to a diffeological space $X$ the set $X$ endowed with the diffeology generated by the set of plots of $X$ of dimension $\le d$. We can easily observe the following:
\begin{itemize}
	\item The functor $\cdot_d$ preserves final structures (see \cite[p. 90]{CSW}).
	\item The functor $\cdot_d$ preserves $\dcal$-embeddings (i.e., if $i: A \longhookrightarrow X$ is a $\dcal$-embedding, then so is $id: A_d \longhookrightarrow X_d$).
	\item $X\cong Z_d$ for some $Z\in \dcal \Leftrightarrow X = X_d \Leftrightarrow {\rm dim}\ X \le d$ (see \cite[pp. 47-48]{IZ}).
\end{itemize}

Next, we define the $d$-coskeleton functor
\[
\cdot_{(d)}: \dcal \longrightarrow \dcal
\]
to assign to a diffeological space $X$ the set $X$ endowed with the diffeology 
\[
D_{X_{(d)}} = \{\text{parametrizations $U \overset{p}{\longrightarrow} X$ such that $U_d \overset{p}{\longrightarrow} X$ is smooth}\}. 
\]
(We can easily check that $D_{X_{(d)}}$ is actually a diffeology on the set $X$.)

We have the natural smooth maps
\[
X_d \overset{id}{\longrightarrow} X \overset{id}{\longrightarrow} X_{(d)}. 
\]
Further, we have the following result.

\begin{prop}\label{adjoint}
	The pair
	\[
	\cdot_d : \dcal \rightleftarrows \dcal : \cdot_{(d)}
	\]
	is an adjoint pair.
	\begin{proof}
		Noting that the diffeology of $A_d$ is final for $\{U_d \xrightarrow{\ p_d\ } A_d \ |\ p\in D_A\}$, we can easily see the result.
	\end{proof}
\end{prop}

A diffeological space $X$ (or its diffeology) is called $d$-{\sl coskeletal} if $X\cong Z_{(d)}$ for some $Z\in \dcal$. $X$ is $d$-coskeletal if and only if $X=X_{(d)}$. Further, we have the following.

\begin{cor}\label{d-coskeletal}
	For a diffeological space $X$, the following are equivalent.
	\begin{itemize}
		\item[{\rm (i)}] $X$ is $d$-coskeletal (resp. $1$-coskeletal).
		\item[{\rm (ii)}] For any $A\in \dcal$ and any set-theoretic map $f:A \longrightarrow X$, $f$ is smooth if and only if $f$ preserves plots of dimensions $\le d$ (resp. smooth curves).
	\end{itemize}
	\begin{proof}
Using the natural bijection $\dcal(A, X_{(d)})\cong \dcal(A_d, X)$, we have the equivalences
		\begin{eqnarray*}
		\mathrm{(ii)} &\Leftrightarrow& \dcal(A, X)\cong \dcal(A_d, X) \ \text{for any } A \in \dcal \\
		&\Leftrightarrow& \dcal(A, X)\cong \dcal(A, X_{(d)}) \ \text{for any } A \in \dcal \\
		\hspace{10em}&\Leftrightarrow& X \cong X_{(d)} \Leftrightarrow \mathrm{(i)}.\hspace{16.5em} \qedhere
		\end{eqnarray*}
	\end{proof}
\end{cor}

\begin{cor}\label{initial}
	The $d$-skeleton functor $\cdot_d: \dcal \longrightarrow \dcal$ and the $d$-coskeleton functor $\cdot_{(d)}: \dcal \longrightarrow \dcal$ preserve final structures and initial structures, respectively, with respect to the underlying set functor.
	\begin{proof}
		Note that the functors $\cdot_d$ and $\cdot_{(d)}$ preserve underlying sets. Then, the result follows from \cite[Proposition 8.7.4]{FK}.
	\end{proof}
\end{cor}
\begin{rem}\label{(co)skeleton}
	The $d$-skeleton and $d$-coskeleton functors appear in the same context as the $d$-skeleton and $d$-coskeleton functors appear in the theory of simplicial sets. More precisely, let $\mathsf{D}$ be the concrete site of Euclidean domains (see \cite[Lemma 4.14]{BH}) and $\mathsf{D}_{\le d}$ be the concrete subsite of $\mathsf{D}$ consisting of Euclidean domains of dimension $\le d$. Then, the $d$-truncation functor 
	\[
	\mathsf{D}_{\le d} Sp \overset{{\rm tr}_d}{\longleftarrow} \mathsf{D}Sp = \dcal 
	\]
	is defined to be precomposition by the inclusion functor $i_d: \mathsf{D}_{\le d} \longhookrightarrow \mathsf{D}$, where $\mathsf{C}Sp$ denotes the category of concrete sheaves on a concrete site $\mathsf{C}$ (or equivalently, the category of quasi-spaces on a concrete site $\mathsf{C}$) (see \cite[Section 4]{BH}, \cite{Dubuc}). We can construct the left adjoint ${\rm sk}_d$ and the right adjoint ${\rm cosk}_d$ of ${\rm tr}_d$ in a similar way to the functors $\cdot_d$ and $\cdot_{(d)}$:
	\begin{equation*}
		\begin{tikzcd}
			\mathsf{D}_{\le d} Sp \arrow[r,bend right=60,"\perp" yshift=10pt ,"{\rm cosk}_d"swap] \arrow[r,bend left=60, "{\rm sk}_d", "\perp"'  yshift=-5pt] & \mathsf{D}Sp = \dcal \arrow[l, "{\rm tr}_d"swap]
		\end{tikzcd}
	\end{equation*}
	Then, we have
	\begin{align*}
		&\cdot_d = {\rm sk}_d \circ {\rm tr}_d : \dcal \longrightarrow \dcal,\\
		&\cdot_{(d)} = {\rm cosk}_d \circ {\rm tr}_d : \dcal \longrightarrow \dcal,
	\end{align*}
	which imply Proposition \ref{adjoint} via \cite[Exercise 4.1.iii]{Riehl17}. (Cf. \cite[Example 6.2.12]{Riehl17}).
\end{rem}

Let $X$ be a diffeological space. Then, $X_0$ is just the set $X$ endowed with discrete diffeology, and hence $X_{(0)}$ is just the set $X$ endowed with indiscrete diffeology. On the other hand, $X_\infty$ is just the diffeological space $X$ itself, and hence so is $X_{(\infty)}$. Thus, from now on, {\bf we assume that} $0<d<\infty$.

The following result is basic even for the study of skeletal diffeological spaces.
\begin{prop}\label{mfdcosk}
	Every $C^\infty$-manifold $M$ is $1$-coskeletal, and hence $d$-coskeletal for any $d>0$. 
	\begin{proof}
		Recall from \cite[27.2]{KM} that smooth maps between $C^\infty$-manifolds are just set-theoretic maps preserving smooth curves. Then, the result is obvious.
	\end{proof}
\end{prop}

\begin{proof}[Proof of Proposition \ref{k-trun}]
	{\rm (1)} The map $id:X_{d}\longrightarrow X$ is obviously smooth. Since the usual topology of an open set of $\rbb^n$ is just the final topology for the smooth curves \cite[Lemma 2.12]{origin}, $\widetilde{id}: \widetilde{X_d} \longrightarrow \widetilde{X}$ is a homeomorphism.
	
	{\rm (2)} By Propositions \ref{adjoint} and \ref{mfdcosk}, we have the (set-theoretic) equalities
	\[
	\dcal(X_d, N) = \dcal(X, N_{(d)}) = \dcal(X,N).
	\]
	We can see from Corollary \ref{d-coskeletal} that the bijection $\dcal(X_d, N) \overset{id^\sharp}{\longleftarrow} \dcal(X,N)$ is a diffeomorphism.
	
	{\rm (3)} By Parts 1 and 2, $\widetilde{X_d} = \widetilde{X}$ and $\dcal(X_d, \rbb) = \dcal(X,\rbb)$ hold. Hence, the result is obvious.\qedhere
\end{proof}
\begin{rem}\label{(co)unit}
	For $X_{(d)}$ and the map $id:X\longrightarrow X_{(d)}$, the obvious analogue of Proposition \ref{k-trun} holds; the proof is similar.
\end{rem}

\begin{exa}\label{(co)skeletal}
	For a $C^\infty$-manifold $M$, the following equivalences hold:
	\begin{align*}
		{\rm dim}\ M\le d & \Leftrightarrow M \ {\rm is }\ d\text{-skeletal (i.e., } M_d=M)\\
		& \Leftrightarrow M_d \ {\rm is }\ d\text{-coskeletal (i.e., } M_d=(M_d)_{(d)}).
	\end{align*}
\end{exa}

We end this section with the following remark.

\begin{rem}\label{final}
	We can observe that the $d$-skeleton functor $\cdot_d$ does not preserve initial structures. In fact, the ordinary diffeology of $\rbb^{d+1} (=\rbb^d \times \rbb)$ is initial for the projections $p_{\rbb^d}: \rbb^{d+1}\longrightarrow \rbb^d$ and $p_\rbb: \rbb^{d+1} \longrightarrow \rbb$. However, the diffeology of $(\rbb^{d+1})_d$ is not initial for the projections $(p_{\rbb^d})_d$ and $(p_\rbb)_d$. This example also shows that $\cdot_d$ does not preserve finite products.
\end{rem}

\section{Proofs of Theorem \ref{immersion}, Corollaries \ref{auto}-\ref{loccontr}, and Theorem \ref{notVD}}
In this section, we prove Theorem \ref{immersion}, Corollaries \ref{auto}-\ref{loccontr}, and Theorem \ref{notVD}.

As mentioned in Section 1, Theorem \ref{immersion} and its corollaries show that smooth maps between $C^{\infty}$-manifolds behave well for the skeleton functor $\cdot_{d}:\mathcal{D}\longrightarrow\mathcal{D}$, but that smooth homotopies between them don't. This is due to the fact that the $d$-skeleton functor $\cdot_d$ does not preserve finite products (Remark \ref{final}).

\begin{proof}[Proof of Theorem \ref{immersion}]
	{\rm (1)} Since $C^\infty$ is a full subcategory of $\dcal$, $C^\infty(M,N) = \dcal(M,N)$ holds. Since $\cdot_d: \dcal \longrightarrow \dcal$ is compatible with the underlying set functor $\dcal \longrightarrow Set$, it is faithful. Thus, we have the natural maps 
	\[
	C^\infty (M,N) \underset{=}{\longrightarrow} \dcal(M,N) \longhookrightarrow \dcal(M_d, N_d).
	\]
	Recall from \cite[27.2]{KM} that smooth maps between $M$ and $N$ are just set-theoretic maps preserving smooth curves. Then, we also see that $C^\infty(M,N) = \dcal(M_d, N_d)$.
	
	{\rm (2)} A map $\iota:D^{d+1}\longrightarrow M$ is called an embedding of $D^{d+1}$ into $M$ if it extends to an ordinary embedding of some open neighborhood $U$ of $D^{d+1}$ in $\rbb^{d+1}$ into $M$. Given an embedding $\iota:D^{d+1}\longrightarrow M,\ K:=\iota(D^{d+1})$ and $\mathring{K}:=\iota(\mathring{D}^{d+1})$ are called the smoothly embedded $(d+1)$-dimensional disk in $M$ and its interior, respectively.
	
	We must show that if $F:M_{d}\times\mathbb{R}\longrightarrow N_{d}$ is a smooth homotopy with $F_0=f$, then $F_1=f$. (Here, $F_t:= F(\cdot, t)$.) Note that since $F$ maps smooth curves of $M\times\mathbb{R}$ to smooth curves of $N$, $F:M\times \rbb\longrightarrow N$ is also a smooth homotopy with respect to the ordinary diffeologies of $M\times\mathbb{R}$ and $N$. Thus, we show that for any smoothly embedded $(d+1)$-dimensional disk $K$ of $M$, the composite
	\[
	K \times [0,1] \longhookrightarrow M \times \rbb \overset{F}{\longrightarrow} N 
	\]
	is a constant smooth homotopy.
	
	Consider the nonempty set
	\[
	S:=\{s\in[0,1]\hspace{0.125cm}|\hspace{0.125cm} F_t |_K =f|_{K}\ {\rm for}\ 0 \le t\leq s\}\subset[0,1] 
	\]
	and set $s_0=\rm{sup}\ \it{S}$. Suppose that $s_0<1$. Then, since $F_t|_K = f|_K$ for $0\le t \le s_0$, there exists $\epsilon>0$ such that $T_x(F_t|_{\mathring{K}})$ has the maximal rank for $x\in \mathring{K}$ and $t\in[0,s_0 \text{+}\epsilon]$ (recall the definition of a smoothly embedded ($d+1$)-dimensional disk). Since $F(x_{0},\cdot):[0 ,s_0\text{+}\epsilon]\longrightarrow N$ is nonconstant for some $x_{0}\in \mathring{K}$, there exists a (linearly) embedded $d$-dimensional open disk $K'$ in $\mathring{K}$ such that
	\[
	F':=F|_{K' \times (0, s_0 + \epsilon)}: K'\times(0 ,s_0\text{+}\epsilon) \longrightarrow N
	\]
	satisfies rank $T_{(x_0, t_0)}F'=d+1$ for some $t_{0}\in(0,s_0\text{+}\epsilon)$. Noting that the inclusion $K' \times (0, s_0+\epsilon) \longhookrightarrow M_d \times \rbb$ is smooth, we see that $K' \times (0, s_0 + \epsilon) \overset{F'}{\longrightarrow} N_d$ is also smooth, which is a contradiction.\qedhere
\end{proof}

A point $x$ of a diffeological space $X$ is called {\sl isolated} if  $\{x\}$ is a connected component of $X$ (see \cite[5.7]{IZ}). Theorem \ref{immersion}(2) states that every immersion $f:M\longrightarrow N$ is isolated in $\dcal(M_d, N_d)$ for $d<{\rm dim}\ M$. Corollaries \ref{auto}-\ref{loccontr} are derived from the isolatedness of self-diffeomorphisms (especially, identities).

\begin{proof}[Proof of Corollary \ref{auto}]
	We see from Theorem \ref{immersion}(1) that ${\rm Aut}_\dcal(M_d)$ is the (set-theoretic) group of diffeomorphisms of $M$. Since the inclusion Aut$_{\mathcal{D}}(M_d)\hookrightarrow\mathcal{D}(M_d,M_d)$ is smooth, ${\rm Aut}_\dcal(M_d)$ has the discrete diffeology by Theorem \ref{immersion}(2). 
\end{proof}
\begin{proof}[Proof of Corollary \ref{homotopydiffeo}]
	(i) $\Leftrightarrow$ (ii) For smooth maps $f:M_d \longrightarrow N_d$ and $g:N_d\longrightarrow M_d$, the equivalence
	\[
	g\circ f \simeq_\dcal 1_{M_d},\ f\circ g \simeq_\dcal 1_{N_d} \Leftrightarrow g\circ f = 1_{M_d},\ f\circ g = 1_{N_d}
	\]
	holds by Theorem \ref{immersion}(2). From this, the equivalence at issue follows.
	
	(ii) $\Leftrightarrow$ (iii) The equivalence follows from Theorem \ref{immersion}(1). \qedhere
\end{proof}

\begin{proof}[Proof of Corollary \ref{loccontr}]
	Note that the functor $\cdot_d$ preserves open $\dcal$-embeddings and that $1_N$ is isolated in $\dcal(N_d,N_d)$ for any $C^\infty$-manifold $N$ of dimension $>d$. Then, the result is obvious.
\end{proof}

\begin{rem}\label{mfd}
	{\rm (1)} Let $M$ and $N$ be second countable finite-dimensional $C^\infty$-manifolds. Recall the exponential law $\dcal(\rbb, \dcal(M,N)) \cong \dcal(\rbb \times M, N)$ and the integrability theorem for vector fields \cite[8.10]{BJ}. Then, we can easily see that neither Theorem \ref{immersion}(2) nor Corollary \ref{auto} holds for $C^\infty$-manifolds with ordinary diffeology. In fact, the set $\mathfrak{C}^\infty(M,N)$ of smooth maps from $M$ to $N$ admits a canonical infinite-dimensional $C^\infty$-manifold structure, for which the map $id:\mathfrak{C}^\infty(M,N)\longrightarrow \dcal(M,N)$ is smooth. Further, if $M$ is compact, then $id:\mathfrak{C}^\infty(M,N)\longrightarrow \dcal(M,N)$ is a diffeomorphism. (See \cite[Theorems 42.1 and 42.14]{KM}).
	
	Since the group ${\rm Diff}(M)$ of diffeomorphisms of $M$ is an open set of $\mathfrak{C}^\infty(M,M)$, it is also an infinite-dimensional $C^\infty$-manifold (actually an infinite-dimensional Lie group) (see \cite[Theorem 43.1]{KM}). Since ${\rm Aut}_\dcal(M)$ is a diffeological subspace of $\dcal(M,M)$ by the inverse function theorem, the map $id: {\rm Diff}(M) \longrightarrow {\rm Aut}_\dcal(M)$ is smooth. Further, if $M$ is compact, then  $id: {\rm Diff}(M) \longrightarrow {\rm Aut}_\dcal(M)$ is a diffeomorphism.
	
	{\rm (2)} A smooth homotopy type of $C^\infty$-manifolds (with ordinary diffeology) can contain many diffeomorphism types (cf. Corollary \ref{homotopydiffeo}).
	
	{\rm (3)} Since $\dcal$ is cartesian closed (Proposition \ref{conven}(2)), $\dcal$ itself is a $\dcal$-category. Theorem \ref{immersion}, along with Part 1 shows that $\dcal(X,Y) \overset{\cdot_d}{\longrightarrow} \dcal(X_d, Y_d)$ need not be smooth, and hence that the functor $\cdot_d: \dcal \longrightarrow \dcal$ cannot be enriched over $\dcal$.
\end{rem}

\begin{rem}\label{trunhomotopy}
	Theorem \ref{immersion}(1) can be strengthened as follows. Suppose that $l \ge 0$ and $d+ l \le e$. Then, $C^\infty(M,N) = \dcal(M,N) = \dcal(M_d,N_e)$ hold in $Set$. Further, $id:\dcal(M_d, N_e) \longrightarrow \dcal(M,N)$ is a smooth bijection, which yields the diffeomorphism $id: \dcal(M_d, N_e)_l \longrightarrow \dcal(M, N)_l$ between diffeological spaces with $l$-skeletal diffeologies.
	
	\begin{proof}[Proof of Theorem \ref{notVD}]
		{\rm (1)} The result follows from \cite[Corollary 11.3]{smh} and \cite[Theorem 3.2]{CW}.
		
		{\rm (2)} We first deal with the case where condition (i) is satisfied, and then use the argument there to deal with the case where condition (ii) is satisfied.\\ 
		{\sl The case where condition (i) is satisfied.} Let $M$ be a $c^\infty$-open set $V$ of a convenient vector space $E$. By translation, we may assume that $x_0 = 0 \in V$. Choose a smoothly contractible open neighborhood $V'$ of $0$ in $V$ (see \cite[Lemma 4.17]{KM}).
		
		We choose a closed linear subspace $F$ of $E$ and an isomorphism $E\cong \rbb^{d+1} \oplus F$ (see \cite[p. 50]{Sch}), and define the smooth map $p: V \longrightarrow \rbb^{d+1}$ to be the composite
		\[
		V \longhookrightarrow E \cong \rbb^{d+1} \oplus F \overset{proj}{\longrightarrow} \rbb^{d+1}.
		\]
		For sufficiently small $r>0$, we define the pointed smooth map $f_r:(S^d, (1,0,\ldots,0))\longrightarrow (V_d, 0)$ by the following conditions:
		\begin{itemize}
			\item $f_r$ corestricts to $V'$.
			\item The composite
			\[
			S^d \overset{f_r}{\longrightarrow} V_d \overset{id}{\longrightarrow} V \longhookrightarrow E \cong \rbb^{d+1} \oplus F
			\]
			is given  by
			\[
			x \longmapsto (r\cdot (x-(1,0,\ldots,0)), 0). 
			\]
		\end{itemize}
		By the construction, $[\widetilde{f_r}]= 0$ in $\pi_d(\widetilde{V},0)$ for any $r$. Suppose that $[f_r]=[f_{r'}]$ in $\pi^\dcal_d (V_d,0)$ for some $r$, $r'$ with $r<r'$. Then, there exists a smooth homotopy 
		\[
		S^d \times \rbb \overset{\widehat{f}}{\longrightarrow} V_d 
		\]
		such that
		\[
		\widehat{f}(\cdot, 0) = f_r, \widehat{f}(\cdot, 1)=f_{r'}, \text{ and } \widehat{f}((1,0,\ldots,0), \cdot)=0.
		\]
		Observe that the image of the composite
		\[
		S^d \times \rbb \overset{\hat{f}}{\longrightarrow} V_d \overset{id}{\longrightarrow} V \overset{p}{\longrightarrow} \rbb^{d+1} 
		\]
		must contain $D^{d+1}_{r'}(-r', 0, \ldots, 0)-D^{d+1}_r(-r,0,\ldots,0)$, where $D^{r+1}_r(-r,0,\ldots,0)$ is the $(d+1)$-dimensional disk with center $(-r, 0, \ldots, 0)$ and radius $r$. Since the composite factors through $V_d$, this is a contradiction by Sard's theorem.\\
		{\sl The case where condition {\rm (ii)} is satisfied}. Choose a chart of $M$ centered at $x_0$ 
		\[
		M \supset U \overset{u}{\underset{\cong}{\longrightarrow}} V \subset E 
		\]
		and use it to identify $U$ with $V$. Then, we can regard the smooth map $p:V\longrightarrow \rbb^{d+1}$ constructed above as a smooth map defined on $U$.
		
		Next, by the $C^\infty$-regularity of $M$, we choose a smooth function $\phi:M \longrightarrow [0,1]$ such that $\phi(x_0)=1$ and supp $\phi \subset U$. By postcomposing an appropriate smooth map $[0,1] \longrightarrow [0,1]$, we may assume that $\phi \equiv 1$ on some open neighborhood $U'$  of $x_0$ in $U$. Further, we may assume that $U' = u^{-1} V'$ by making $V'$ smaller if necessary. 
		
		Then, we can define the smooth map $q:M\longrightarrow \rbb^{d+1}$ by
		\begin{equation*}
			q(x) = 
			\begin{cases}
				\phi(x)\cdot p(x) & \text{for } x\in U,\\
				0 & \text{for } x \notin U.
			\end{cases}
		\end{equation*}
		For sufficiently small $r>0$, we define the pointed smooth map $f_r:(S^d,(1,0,\ldots,0))\longrightarrow (M_d, x_0)$ as above. Then, it is obvious that $[\widetilde{f}_r] = 0$ is $\pi_d(\widetilde{M}, x_0)$ for any $r$. By using $q$ instead of $p$, we can see that if $r\neq r'$, then $[f_r]\neq [f_{r'}]$ in $\pi^\dcal_d (M_d, x_0)$.
	\end{proof}
\end{rem}

\begin{rem}\label{homology/homotopy}
By the argument in the proof of Theorem \ref{notVD}, we can show that under the same assumption, the following hold:
\begin{itemize}
\item The kernel of the homomorphism
$$
H_{d}(M_{d};A) \longrightarrow H_{d}(\widetilde{M};A)
$$
is uncountable for any coefficient group $A\ (\neq 0)$ (see \cite[Section 3.1]{smh}).
\item The kernel of the homomorphism
$$
\pi^{\dcal}_{i} (M_d, x_0) \longrightarrow \pi_i(\widetilde{M}, x_0)
$$
is uncountable for any $i$ such that $\pi_{i}(S^{d}) \neq 0$.
\end{itemize}
\end{rem}

We end this section by describing where $C^{\infty}$-manifolds with skeletal diffeology go in Venn diagram (A.1) of \cite[Appendix A]{smh}.

\begin{rem}\label{WV}
	Recall the subclasses $\wcal_\dcal$, $\vcal_\dcal$, and $\widetilde{}\wcal_\czero$ of $\dcal$ from \cite[Section 1.3 and Appendix A]{smh}. Then, Theorem \ref{notVD} implies the following: If $M$ is a $c^\infty$-open set of a convenient vector space or a $C^\infty$-regular $C^\infty$-manifold,  then $M_d$ is not in $\vcal_\dcal$ for $0<d<{\rm dim}\ M$.
	
	Let $M$ be a paracompact finite-dimensional $C^\infty$-manifold. Since $M$ is hereditarily $C^\infty$-paracompact (see Remark \ref{Cparacpt}),
	\[
	M \in \wcal_\dcal \subset \vcal_\dcal \cap \widetilde{}\wcal_\czero
	\]
	holds by \cite[Theorem 11.1 and Corollary 1.6]{smh}. On the other hand, since $M$ is $C^\infty$-regular,
	\[
	M_d \in \widetilde{} \wcal_\czero \backslash \vcal_\dcal
	\]
	holds for $0<d<{\rm dim}\ M$.
\end{rem}




\renewcommand{\thesection}{Appendix}
\section{Coskeletal diffeological spaces}
\renewcommand{\thesection}{A}
In this appendix, we establish some results on coskeletal diffeological spaces.

The $d$-coskeleton functor $\cdot_{(d)}: \dcal \longrightarrow \dcal$ preserves limits and initial structures (Proposition \ref{adjoint} and Corollary \ref{initial}), and the obvious analogue of Proposition \ref{k-trun} holds for $X_{(d)}$ and $id : X \longrightarrow X_{(d)}$ (see Remark \ref{(co)unit}). Thus, the functor $\cdot_{(d)}: \dcal \longrightarrow \dcal$ preserves group objects, locally trivial (principal) bundles, and $\dcal$-numerable (principal) bundles (see \cite[Definitions 2.6 and 2.13]{CW21}).

Further, we have the following result.

\begin{prop}\label{D(A,X)}
	Let $0<d<\infty$. If $X$ is a $d$-coskeletal diffeological space, then $\dcal(A,X)$ is also $d$-coskeletal for any $A\in \dcal$. 
	\begin{proof}
		Let $p: U \longrightarrow \dcal(A,X)$ be a parametrization of $\dcal(A,X)$, and $\widehat{p}: U \times A \longrightarrow X$ be the map corresponding to $p$ via the (set-theoretic) exponential law.
		
		Since $X$ is $d$-coskeletal, we can use Corollary \ref{d-coskeletal} to see that the following conditions are equivalent:
		\begin{itemize}
			\item[{\rm (i)}] $p:U\longrightarrow \dcal(A,X)$ is a plot of $\dcal(A,X)$.
			\item[{\rm (ii)}] $\widehat{p}: U \times A \longrightarrow X$ is smooth.
			\item[{\rm (iii)}] The composite $V\xrightarrow{\ (q, r)\ } U \times A \overset{\widehat{p}}{\longrightarrow} X$ is smooth for any $q\in D_U$, $r\in D_A$ with a common domain of dimension $\le d$.
		\end{itemize}
		Further, we can observe that (iii) is equivalent to
		\begin{itemize}
			\item[{\rm (iv)}] The composite $V \times A \xrightarrow{\ q\times 1_A \ } U \times A \overset{\widehat{p}}{\longrightarrow} X$ is smooth for any plot $q\in D_U$ of dimension $\le d$.
		\end{itemize}
		By the cartesian closedness of $\dcal$, (iv) is equivalent to
		\begin{itemize}
			\item[{\rm (v)}] The composite $V \overset{q}{\longrightarrow} U \overset{p}{\longrightarrow} \dcal(A,X)$ is smooth for any plot $q\in D_U$ of dimension $\le d$, 
		\end{itemize}
		and hence to
		\begin{itemize}
			\item[{\rm (vi)}] $p:U\longrightarrow \dcal(A,X)$ is a plot of $\dcal(A,X)_{(d)}$.\qedhere
		\end{itemize}
	\end{proof}
\end{prop}

\begin{cor}\label{1-coskeletal}
	Let $M$ be a $C^\infty$-manifold.
	\begin{itemize}
		\item[{(1)}] $\dcal(A,M)$ is $1$-coskeletal for any $A\in \dcal$.
		\item[{(2)}] ${\rm Aut}_\dcal(M)$ is $1$-coskeletal.
	\end{itemize}
	\begin{proof}
		The result follows from Propositions \ref{mfdcosk} and \ref{D(A,X)} and Corollary \ref{initial} (recall the definition of the automorphism group ${\rm Aut}_\dcal(F)$ from Section 1).
	\end{proof}
\end{cor}

By Corollaries \ref{1-coskeletal} and \ref{d-coskeletal}, a set-theoretic map $f:S\longrightarrow \dcal(A,M)$ is smooth if and only if $f$ preserves smooth curves.

\begin{rem}\label{coskmfd}
	If $A$ is a compact $C^\infty$-manifold and $M$ is a paracompact finite-dimensional $C^\infty$-manifold, then the $1$-coskeletal diffeological space $\dcal(A,M)$ is a $C^\infty$-manifold. Similarly, if $M$ is a compact $C^\infty$-manifold, then the $1$-coskeletal diffeological group ${\rm Aut}_\dcal(M)$ is a Lie group. (See Remark \ref{mfd}(1).)
\end{rem}

\section*{Acknowledgements}
In the conference ``Building-up Differentiable Homotopy Theory 2019" held at Kyusyu University, Professor Iglesias-Zemmour asked me where the sphere with wire diffeology would go in Venn diagram (A.1) in \cite[Appendix A]{smh} during my talk. Answering this question is one of the motives for this work (see Remark \ref{WV}). I would like to thank Professor Iglesias-Zemmour for his insightful qustion.

\renewcommand{\thesection}{A}


\bibliographystyle{elsarticle-num}
\bibliography{<your-bib-database>}

\begin{thebibliography}{00}



	
	
	
	\bibitem{BH}
	J. Baez and A. Hoffnung, {\sl Convenient categories of smooth
		spaces}, Transactions of the American Mathematical Society \textbf{363} (2011), no. 11, 5789-5825.
	
	
	
	
	\bibitem{BJ}
	T. Br\"{o}cker and K. J\"{a}nich, {\sl Introduction to differential topology,} Cambridge University Press, (1982).
	
	
	\bibitem{CSW}
	J. D. Christensen, G. Sinnamon, and E. Wu, {\sl The D-topology
		for diffeological spaces,} Pacific Journal of Mathematics \textbf{272} (2014), no. 1, 87-110.
	
	\bibitem{CW}
	J. D. Christensen and E. Wu, {\sl The homotopy theory of diffeological
		spaces,} New York J. Math \textbf{20} (2014) 1269-1303.
	
	
	
	
	\bibitem{CW16a}
	J. D. Christensen and E. Wu, {\sl Tangent spaces and tangent bundles for diffeological spaces}, Cahiers de Topologie et Geom trie Diff rentielle Cat gorigues 57(1) (2016), 3-50.
	\bibitem{CW20}
	J. D. Christensen and E. Wu, {\sl Exterior bundles in diffeology}, arXiv preprint arXiv:2009.01770 (2020).
	
	\bibitem{CW21}
	J. D. Christensen and E. Wu, {\sl Smooth classifying spaces}, Israel Journal of Mathematics 241.2 (2021): 911-954.
	
	\bibitem{Dubuc}
	E. J. Dubuc, {\sl Concrete quasitopoi}, Applications of sheaves, Springer, Berlin, Heidelberg, (1979) 239-254.
	
	
	
	
	
	
	
	\bibitem{FK}
	A. Fr\"{o}licher and A. Kriegl, {\sl Linear spaces and differentiation theory}, Vol. 13, John Wiley and Sons Inc, (1988).
	
	\bibitem{Gauld}
	D. B. Gauld, {\sl Topological properties of manifolds}, The American Mathematical Monthly 81.6 (1974): 633-636.
	
	
	
	
	
	
	
	
	
	
	
	\bibitem{IZ}
	P. Iglesias-Zemmour, {\sl Diffeology,} Vol. 185, American Mathematical Soc, (2013).
	
	
	
	
	
	
	
	
	
	\bibitem{origin}
	H. Kihara, {\sl Model category of diffeological spaces,} Journal of Homotopy and Related Structures, 14.1 (2019): 51-90.
	
	\bibitem{smh}
	H. Kihara, \textit{Smooth Homotopy of Infinite-Dimensional $ C^{\infty} $-Manifolds}, to appear in Memoirs of the American Mathematical Society, available at arXiv:2002.03618 (2020).
	
	\bibitem{deRham}
	H. Kihara, {\sl De Rham calculus on diffeological spaces, preprint.}
	
	
	\bibitem{KM}
	A. Kriegl and P. W. Michor, {\sl The convenient setting of global
		analysis,} Vol. 53, American Mathematical Society (1997).
	
	
	
	
	
	
	
	
	
	
	
	
	
	
	
	
	
	\bibitem{Riehl17}
	E. Riehl, {\sl Category theory in context}, Courier Dover Publications, 2017.
	
	
	\bibitem{Sch}
	H. H. Schaefer and M. P. Wolff, {\sl Topological vector spaces}, 2nd ed., Graduate Texts in Mathematics, vol. 3, Springer-Verlag, New York, 1999.

\end{thebibliography}



\end{document}